\documentclass{article}
\usepackage{amsmath,amssymb}
\usepackage{graphicx}
\usepackage[font=small,labelfont=bf]{caption}
\usepackage{epstopdf}

\begin{document}

\title{A Note on the Sparing Number of the Sieve Graphs of Certain Graphs\thanks{%
Mathematics Subject Classifications: 05C78.}}
\date{{\small }}
\author{N. K. Sudev\thanks{%
Department of Mathematics, Vidya Academy of Science and Technology, Thalakkottukara, Thrissur-680501, Kerala, India,
}\ , K. A. Germina\thanks{%
PG \& Research Department of Mathematics, Mary Matha Arts \& Science College, Mananthavady, Wayanad-670645, Kerala, India.}}
\maketitle

\begin{abstract}
Let $\mathbb{N}_0$ denote the set of all non-negative integers and $\mathcal{P}(\mathbb{N}_0)$ be its power set. An integer additive set-indexer (IASI) of a given graph $G$ is an injective function $f:V(G)\to \mathcal{P}(\mathbb{N}_0)$ such that the induced function $f^+:E(G) \to \mathcal{P}(\mathbb{N}_0)$ defined by $f^+ (uv) = f(u)+ f(v)$ is also injective. An IASI $f$ of a graph $G$ is said to be a weak IASI of $G$ if $|f^+(uv)|=\max(|f(u)|,|f(v)|)$ for all $u,v\in V(G)$. A graph which admits a weak IASI may be called a weak IASI graph. The sparing number of a graph $G$ is the minimum number of edges with singleton set-labels required for a graph $G$ to admit a weak IASI.  In this paper, we introduce the notion of $k$-sieve graphs of a given graph and study their sparing numbers.
\end{abstract}

\section{Introduction}

For all  terms and definitions, not defined specifically in this paper, we refer to \cite{FH} and \cite{DBW} and for different graph classes, we refer to \cite{BLS}. Unless mentioned otherwise, all graphs considered here are simple, finite and have no isolated vertices.

The notion of a set-valued graph has been introduced in \cite{A1} as a graph, the labels of whose vertices and edges are the subsets of a given set. Since then, several studies have been done on set-valuations of graphs. The {\em sumset} of two non-empty sets $A, B$, denoted by  $A+B$, is defined as $A + B = \{a+b: a \in A, b \in B\}$. Using the terminology of sumsets of sets, the notion of an integer additive set-indexer of a given graph is introduced in \cite{GA} as follows. 

Let $\mathbb{N}_0$ denote the set of all non-negative integers and $\mathcal{P}(\mathbb{N}_0)$ be its power set. An {\em integer additive set-indexer} (IASI, in short) of a graph $G$ is an injective function $f:V(G)\to \mathcal{P}(\mathbb{N}_0)$ such that the induced function $f^+:E(G) \to \mathcal{P}(\mathbb{N}_0)$ defined by $f^+ (uv) = f(u)+ f(v)$ is also injective. 

The cardinality of the set-label of an element (vertex or edge) of a graph $G$ is called the {\em set-indexing number} of that element.

\smallskip

LEMMA 1.1. \cite{GS1} Let $A$ and $B$ be two non-empty finite sets of non-negative integers. Then, $\max(|A|,|B|) \le |A+B|\le |A|\,|B|$. Therefore, for any integer additive set-indexer $f$ of a graph $G$, we have $\max(|f(u)|,\, |f(v)|)\le |f^+(uv)|= |f(u)+f(v)| \le |f(u)|\, |f(v)|$, where $uv\in E(G)$.

\smallskip

DEFINITION 1.2. \cite{GS1} An IASI $f$ of a graph $G$ is said to be a {\em weak IASI} if $|f^+(uv)|=|f(u)+f(v)|=\max(|f(u)|,\,|f(v)|)$ for all $u,v\in V(G)$. A graph which admits a weak IASI is called a {\em weak IASI graph}. A weak  IASI $f$ is said to be {\em weakly $k$-uniform IASI} if $|f^+(uv)|=k$, for all $u,v\in V(G)$ and for some positive integer $k$.

\smallskip

If $A$ and $B$ are two non-empty sets of non-negative integers, then $|A+B|=|A|$ if and only if $|B|=1$ and $|A+B|=|B|$ if and only if $|A|=1$. Hence, we have

\smallskip

THEOREM 1.3. \cite{GS1}  A graph $G$ admits a weak IASI if and only if at least one end vertex of every edge of $G$ has a singleton set-label.
 
\smallskip

DEFINITION 1.4. \cite{GS3} A {\em mono-indexed element} (a vertex or an edge) of an IASI graph $G$ is an element of $G$ whose set-indexing number is $1$. The {\em sparing number} of a graph $G$ is defined to be the minimum number of mono-indexed edges required for $G$ to admit a weak IASI and is denoted by $\varphi(G)$.

\smallskip

THEOREM 1.5. \cite{GS3} An odd cycle $C_n$ contains odd number of mono-indexed edges and an even cycle contains an even number of mono-indexed edges.

\smallskip

THEOREM 1.6. \cite{GS3} The sparing number of an odd cycle $C_n$ is $1$ and that of an even cycle is $0$.

\smallskip

THEOREM 1.7. \cite{GS3} The sparing number of a bipartite graph is $0$.

\smallskip

THEOREM 1.8. \cite{GS3} The sparing number of a complete graph $K_n$ is $\frac{1}{2}(n-1)(n-2)$.

\smallskip

\noindent Now, recall the definition of graph powers.

\smallskip

DEFINITION 1.9. \cite{BM1} The $r$-th power of a simple graph $G$ is the graph $G^r$ whose vertex set is $V$, two distinct vertices being adjacent in $G^r$ if and only if their distance in $G$ is at most $r$. The graph $ G^2 $ is referred to as the {\em square} of $G$, the graph $ G^3 $ as the {\em cube} of G.

\smallskip

\noindent  The following is an important theorem on graph powers.

\smallskip

THEOREM 1.10. \cite{EWW} If $d$ is the diameter of a graph $G$, then $G^d$ is a complete graph.

\smallskip

\noindent The admissibility of weak IASIs by certain graph classes and graph powers and the determination of their corresponding  sparing numbers have been done in \cite{GS5}, \cite{GS6} and \cite{GS7}. The admissibility of weak IASIs by the graph operations and certain graphs associated with the given IASI graphs have been discussed in \cite{GS0} and \cite{GS4}. As a continuation to these studies, in this paper, we discuss the sparing number of a particular type of graphs obtained by adding some edges to the given graphs according to certain rules.

\section{Sparing Number of the $k$-Sieve of a Graph}

Motivated by the terminology of graph powers, we introduce the notion of a $k$-sieve of a given graph as follows.

\smallskip
 
DEFINITION 2.1. A {\em $k$-sieve graph} or simply a {\em $k$-sieve} of a given graph $G$, denoted by $G^{(k)}$, is the graph obtained by joining the non-adjacent vertices of $G$ which are at a distance $k$. A cycle obtained by joining two vertices of $G$, which are at a distance $k$ in $G$, is called a {\em $k$-ringlet} of the graph $G$.

\smallskip

Note that every $k$-ringlet of a graph is a cycle of length $k+1$. The number of $k$-ringlets in a graph $G$ is the number of distinct $k$-paths in $G$. The number of edges in $G^{(k)}$ that are not in $G$ is the number $k$-ringlets in $G$.

\smallskip

REMARK 2.2.  Note that $G^{(2)} \cong G^2$, the square of the graph $G$. But, for $k>2$, $G^{(k)}$ and $G^k$ are non-isomorphic graphs. The sparing number of the square of certain graphs are studied in \cite{GS7}. Hence, in this paper, we need to consider $k\ge 3$.

\smallskip

PROPOSITION 2.3. Let $l$  be the length of a maximal path in $G$. If $k>l$, then the sparing number of the $k$-sieve of $G$ is equal to the sparing number of $G$ itself. 

\smallskip

PROOF. Given that $l$ is the length of a maximal path in $G$. Hence, for any pair of vertices $x,y$ in $G$, $d(x,y)\le l$. That is, there exists no vertex in $G$ which is at a distance $k$ from another vertex of $G$. Therefore, if $k>l$, then $G^{(k)}\cong G$. Hence, $\varphi(G^{(k)})=\varphi(G)$. 

\smallskip

Invoking the above result, we need to consider the integral values between $3$ and $l$, including both, for $k$. If $k=l$, then the longest path of $G$ becomes a cycle of length $l+1$ in $G^{(k)}$.

\smallskip

We now proceed to determine the sparing number of the sieve graphs of certain other standard graphs. Let us begin with the path graphs. 

\smallskip

THEOREM 2.4. Let $P_n$ be a path of length $n$. Then, for odd integers $k;~2<k\le n$, the sparing number of $P_n^{(k)}$ is $0$ and for even integer $k; ~2<k\le n$,  the sparing number of $P_n^{(k)}$ is 
\begin{equation*}
\varphi(P_n^{(k)})=
\begin{cases}
2(lk+r)-3 & ~\text{if $(lk+r)k+s=n$ where $s\le r-2$ and $r\ge 2$}\\
2(lk+r)-2 & ~\text{if $(lk+r)k+s=n$ where $s=r-1$ and $r\ge 1$}\\
2(lk+r)-1 & ~\text{if $(lk+r)k+s=n$ where $s=r$ and $r\ge 0$}.
\end{cases}
\end{equation*}
where $l,k$ and $s$ are non-negative integers.

\smallskip

PROOF. Let $P_n$ be a path on $n+1$ vertices. Let $V=\{v_1,v_2,v_3,\ldots, v_n, v_{n+1}\}$ be the vertex set of $P_n$. The proof is developed considering various possible  cases as below.

\smallskip

\noindent {\em Case-1:} Let $n$ be an odd integer. Label the vertices of $P_n$ alternately by distinct singleton sets and distinct non-singleton sets. Then, each vertex $v_i$ is adjacent to $v_{i-1}, v_{i+1}$ and to the vertex $v_{i+k}$ in $P_n^{(k)}$. Then, we have the following subcases.

\smallskip
	
\noindent {\em Subcase-1.1:} If $k=n$, then by Proposition 2.3, $P_n^{(k)}=C_{n+1}$. Since $n$ is odd, $P_n^{(k)}$ is even cycle. Then, by Theorem 1.8, $P_n^{(k)}$ has no mono-indexed edges. 

\smallskip
	
\noindent {\em Subcase-1.2:} Let $k<n$. Then, under the set-labeling we defined above, no two adjacent vertices simultaneously have singleton set-labels or non-singleton set-labels. Therefore, no edge in $P_n^{(k)}$ has no mono-indexed edges. Therefore, $P_n^{(k)}$ has no mono-indexed edges if $k\le n$ is an odd integer.

\smallskip
	
\noindent {\em Case-2:} Let $n$ be an even integer.

\smallskip
	
\noindent {\em Subcase-2.1:} If $k=n$, then by Proposition 2.3, $P_n^{(k)}=C_{n+1}$. Since $n$ is even, $P_n^{(k)}$ is odd cycle. Then, by Theorem 1.8, $P_n^{(k)}$ has at least one mono-indexed edge. That is, $\varphi(P_n^{(k)})=1$.

\smallskip
	
\noindent {\em Subcase-2.2:} Let $k<n$. Then, there exists two integers $r$ and $s$ such that $rk+s\le n$, where $0\le s< k$ and $r=0,1,2,\ldots$. We can label the vertices in such a way that no two adjacent vertices have non-singleton set-labels in the following way. 
	
Label the vertices $v_1,v_3, v_5,\ldots, v_{k-1}$ of $P_n^{(k)}$ by distinct non-singleton sets and label $v_2,v_4,v_6,\ldots, v_k$ by distinct singleton sets. Since $v_{k+1}$ is adjacent to $v_1$, $v_{k+1}$ can be labeled only by a singleton set which is not used for labeling any one of the preceding vertices. Then, the edge $v_kv_{k+1}$ is a mono-indexed edge. If $n\le 2k$, then the only mono-indexed edge in $P_n^{(k)}$ is $v_kv_{k+1}$. 
	
If $n>2k$, label the vertices $v_1, v_2,\ldots, v_{k+1}$ as mentioned above and proceed labeling $v_{k+2}, v_{k+4},\ldots, v_{2k}$ by distinct non-singleton sets, that are not used for labeling before, and the vertices $v_{k+2}, v_{k+4},\ldots, v_{2k}$ by distinct singleton sets that are not used before for labeling. Then, $v_{k+1}v_{2k+1}$ is a mono-indexed edge. Since, $v_{k+2}$ is adjacent to $v_{2k+2}$ must be mono-indexed. Therefore, $v_{2k+1}v_{2k+2}$ is also a mono-indexed edge. Proceeding in this way, we arrive at the following cases.

\begin{enumerate}
\item Let $s\le r-2$. Then, $r\ge 2$. Now, we can find a path $P':v_kv_{k+1}v_{2k+1}v_{2k+2}\\v_{3k+2}v_{3k+3}\ldots\ldots,v_{(r-1)k+(r-2)}v_{(r-1)k+(r-1)}$, all of whose elements have singleton set-labels, containing all the mono-indexed edges of $P_n^{(k)}$. The length of the path $P'$ is $2r-3$.
\item If $s=r-1$, then $r\ge 1$. Now, there exists a path $P'':v_kv_{k+1}v_{2k+1}v_{2k+2}v_{3k+2}\\v_{3k+3}\ldots\ldots, v_{(r-1)k+(r-1)}v_{rk+(r-1)}$, all of whose edges are mono-indexed, containing all the mono-indexed edges of $P_n^{(k)}$. The length of the path $P''$ is $2r-2$.
\item If $s=r$, then the path $P''':v_kv_{k+1}v_{2k+1}v_{2k+2}v_{3k+2}v_{3k+3}\ldots, v_{rk+(r-1)}v_{(r+1)k}$, all of whose edges are mono-indexed, containing all the mono-indexed edges in $P_n^{(k)}$. Therefore, the length of $P'''$ is $2r-1$. 
\end{enumerate} %

If $r=s=k$, then $rk+s=(k+1)k$ and if $r=lk$ and $s=k$, then $rk+s=(lk+1)k$ and hence we can proceed the labeling procedure in the same manner as explained above. Then, the result follows.

\smallskip

\begin{figure}[h!]
	\centering
	\includegraphics[scale=0.35]{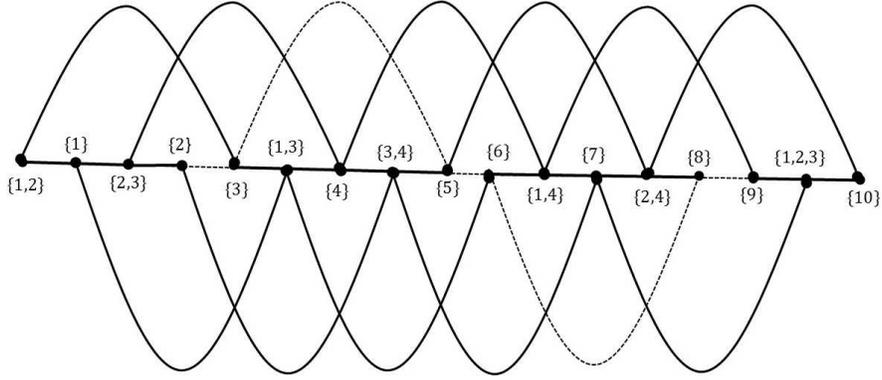}
	\caption{$4$-sieve of a path with a weak IASI defined on it.}\label{G-Pn-k-Circ}
\end{figure}

Figure \ref{G-Pn-k-Circ} illustrates a weak IASI for the $4$-sieve of the path of length $16$. The mono-indexed edges are represented in dotted lines.

\smallskip

The above theorem arouses an interest in determining the sparing number of the $k$-sieve of a cycle. Here, note that a maximal path between any pair of vertices in a cycle $C_n$ is $\lfloor \frac{n}{2} \rfloor$. Therefore, a $k$-sieve exists for $C_n$ if and only if $n\ge 2k+1$. Moreover, $C_{n^(k)}$ is a $4$-regular graph. Then, we have the following results.

\smallskip

THEOREM 2.5. \label{T-Circ-Cyc-E} Let $C_n$ be a Cycle that admits a weak IASI. Then, for an odd integer $k$, $1<k\le l$, 
\begin{eqnarray*}
		\varphi(C_n^{(k)})=
		\begin{cases}
			0 & \text{if $C_n$ is an even cycle}\\
			k+1 & \text{if $C_n$ is an odd cycle}.
		\end{cases}
\end{eqnarray*}
where $l$ is the length of a largest path in $G$.

\smallskip

PROOF. Let $k$ be an odd integer. Then, every $k$-subcycle of $C_n^{(k)}$, obtained by joining the vertices of $C_n$ which are at a distance $k$ in $C_n$, is an even cycle of length $k+1$. Let $C'$ be such an even  cycle of length $k+1$ in $C_n^{(k)}$, which has exactly one edge, say $e'$, which is not in $C_n$. 
	
If $C_n$ is an even cycle, then by Theorem 1.6, it need not contain mono-indexed edges. Therefore, as a result of Theorem 1.5, $e'$ can not be mono-indexed. Therefore, $C'$ does not contain any mono-indexed edges.
	
If $C_n$ is an odd cycle, then by Theorem 1.6, it must have one mono-indexed edge. If $C'$ contains this mono-indexed edge of $C_n$, then as a result of Theorem 1.5, $e'$ must be mono-indexed. There exist such $k$ cycles containing this mono-indexed edge of $C_n$. Therefore, the sparing number of $C_n^{(k)}$ is $k+1$.

\smallskip

\begin{figure}[h!]
	\centering
	\includegraphics[scale=0.4]{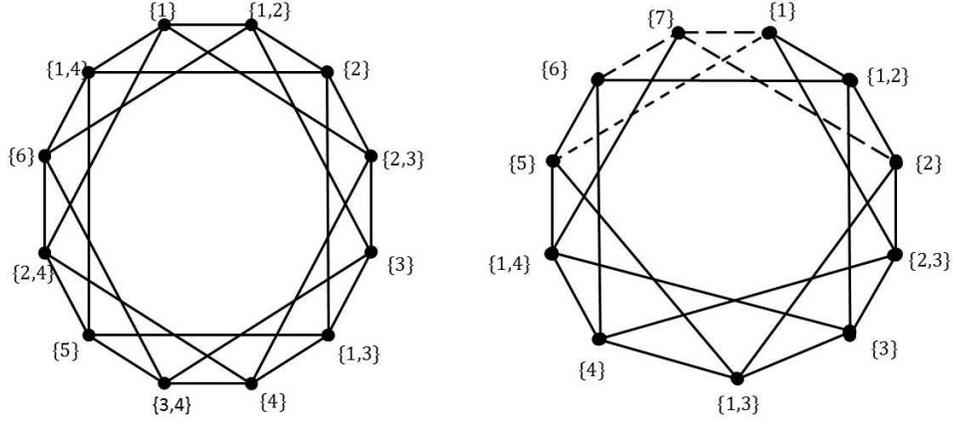}
	\caption{$3$-sieve of $C_{12}$ and a $3$-sieve of $C_{12}$ with weak IASIs defined on them.}\label{G-Circ-k1}
\end{figure}

Figure \ref{G-Circ-k1} illustrates Theorem \ref{T-Circ-Cyc-E}. The first subfigure is the $3$-sieve of an even cycle with a weak IASI defined on it and the second subfigure  is  $3$-sieve of an odd cycle with a weak IASI on it. Mono-indexed edges in the second graph are represented by dotted lines.

\smallskip

\noindent Next let us consider the case when $k$ is an even integer.

\smallskip

THEOREM 2.6 \label{T-Circ-Cyc-O} Let $C_n$ be a cycle of length $n$. For an even integer $k; ~2<k\le n$,  the sparing number of $C_n^{(k)}$ is 
\begin{equation*}
	\varphi(P_n^{(k)})=
	\begin{cases}
	3 & ~\text{if}~~ $n=2k$\\
	2[(lk+r)-2\,\lfloor \frac{(l-1)k+(r-1)}{2} \rfloor] & ~\text{if}~~ n=lk+r\\
	2l & ~\text{if}~~ $n=l(k+1)$\\
\end{cases}
\end{equation*}
where $l,k$ and $r$ are non-negative integers such that $l\ge 2, r<l$.

\smallskip

PROOF. First let $n=2k$. Then , $C_n^{(k)}$ is a cubic graph. Let us begin the labeling process by labeling the first vertex $v_1$ by a non-singleton set and then label the following vertices alternatively by distinct singleton sets and distinct non-singleton sets. Then, the vertex $v_k$ is a mono-indexed vertex. Being adjacent to the vertex $v_1$, $v_{k+1}$ must also be mono-indexed. That is, the edge $v_kv_{k+1}$ is mono-indexed. Now, label the vertex $v_{k+2}$ by a non-singleton set and then label the following vertices alternatively by distinct singleton sets and distinct non-singleton sets. Here, the vertex $v_{2k-1}$ is mono-indexed. Since, $v_{2k}$ is adjacent to the vertex $v_1$, $v_{2k}$ must be mono-indexed. Therefore, the edge $v_{2k-1}v_{2k}$ is mono-indexed. Also, the edge $v_kv_{2k}$ is also mono-indexed. Therefore, the number mono-indexed edges in this case is $3$.  

\smallskip
	
Note that if $n>2k$ the $k$-sieve of every cycle is a $4$-regular graph, for any integer $k$. Then we have the following cases.

\smallskip
	
\noindent {\em Case-1:} Assume that $n=lk+r;~r<l$, $l$ and $r$ being positive integers and $l\ge 2$. Then, the total number of edges in $C_n^{(k)}$ is $|E(C_n^{(k)})|=\frac{1}{2}\sum d(v)=2(lk+r)$. Now, label the vertex $v_1$ by a non-singleton set and then label the remaining vertices by distinct singleton sets and distinct non-singleton sets such that no two adjacent vertices have non-singleton set-labels. Then, the last $k+1$ vertices must be $1$-uniform, as each of them are adjacent to one vertex having a non-singleton set-label. Out of the remaining $(l-1)k+(r-1)$ vertices, $\lfloor \frac{(l-1)k+(r-1)}{2} \rfloor$ vertices have non-singleton set-label. The number of edges that are not mono-indexed is $4\,\lfloor \frac{(l-1)k+(r-1)}{2} \rfloor$. The total number of mono-indexed edges is $2[(lk+r)-2\,\lfloor \frac{(l-1)k+(r-1)}{2} \rfloor]$.

\smallskip
	
\noindent {\em Case-2:} Assume that $n=l(k+1)$, $l$ being a positive integer. Let $\mathfrak{C}$ be a partition of $V(G)$, where each set in $\mathfrak{C}$ contains exactly $k+1$ vertices. Therefore, each set in $\mathfrak{C}$ consists of exactly $\frac{k}{2}$ vertices have non-singleton set-labels and  $1+\frac{k}{2}$ mono-indexed vertices in $C_n^{(k)}$. Therefore, the number of vertices having non-singleton set-labels is $l\frac{k}{2}$. Therefore, the number of edges that are not mono-indexed is $2lk$. The number edges in $C_n^{(k)}$ is $|E(C_n^{(k)})|=\frac{1}{2}\sum d(v)=2l(k+1)$. Therefore, the number of mono-indexed edges in $C_n^{(k)}$ is $2l(k+1)-2lk=2l$.

\smallskip

\begin{figure}[h!]
	\centering
	\includegraphics[scale=0.3]{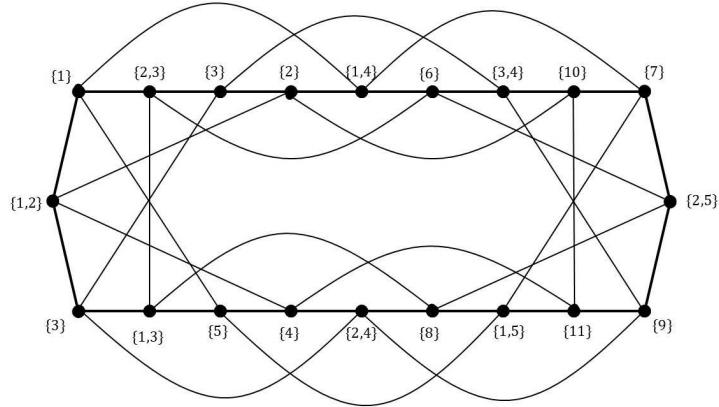}
	\caption{$4$-sieve of $C_{20}$ which is a weak IASI graph.}\label{G-Cyc-k-Circ1}
\end{figure}

Figure \ref{G-Cyc-k-Circ1} illustrates a weak IASI for the $4$-sieve of a cycle on $20$ vertices.

\smallskip

\section{Conclusion}
It can be observed that a complete graph $K_n$ can have a $k$-sieve graph as every vertex of $K_n$ is at a distance $1$ from all other vertices of $G$. Similarly, a complete bipartite graph $K_{m,n}$ (or a complete $r$-partite graph, $K_{n_1,n_2,\ldots,n_r}$, for $r>2$) also does have a $k$-sieve graph, for $k\ge 3$, as any two vertices in $K_{m,n}$ (or in $K_{n_1,n_2,\ldots,n_r}$)are at a distance at most $2$. 

Let $k$ be an odd integer. If $G$ be a tree, then $G^{(k)}$ is a graph all of whose cycles are of length $k+1$, an even integer. Then, $G^{(k)}$ is a bipartite graph. Therefore, by Theorem 1.7, the number of mono-indexed edges in $G^{(k)}$ is also $0$. 

If $G$ be a bipartite graph (containing cycles), then $G$ has no odd cycles. Then, since every $k$-ringlet of $G$ is an even cycle, every cycle in $G^{(k)}$ is of even length. Therefore, $G^{(k)}$ is also a bipartite graph, for odd $k$. Therefore, by Theorem 1.7, the number of mono-indexed edges in $G^{(k)}$ is also $0$.

But, for even integers $k$, to determine the sparing number of $G^{(k)}$, for a graph $G$ which is bipartite (cyclic or acyclic), we need to use Theorem 2.4, Theorem 2.5 and Theorem 2.6, for distinct paths and cycles in $G$.  

Evaluating the sparing number of the $k$-sieves of bipartite graphs, Eulerian graphs, armed crown graphs etc. are some of the open problems in this area. Determining the sparing number of the $k$-sieves of graph operations and graph products are also worth for further exploration.

More properties and characteristics of weak IASIs, both uniform and non-uniform, are yet to be investigated. The problems of establishing the necessary and sufficient conditions for various graphs and graph classes to have certain IASIs still remain to be settled. All these facts highlight a great scope for further studies in this area.

\smallskip

\section*{Acknowledgement}

The authors gratefully acknowledge the critical comments and suggestions the anonymous referee, which helped to improve the quality of the paper.

\end{document}